  \documentstyle[12pt]{article}

\makeatletter

\ifnum\@ptsize=0  \fi
\ifnum\@ptsize=2  \fi
\advance\textheight by 0.7cm
\advance\voffset by -0.3cm
 \newcommand{\qed}{\hspace*{\fill}Q.E.D.\vskip12pt}

\def\Pn#1{{\bf P}^{#1}}

\def\Gr #1#2{{\bf G}(#1,#2)} 

\def\Pd{{\bf P}(\wedge^2V)}
\def\Pr{{\bf P}(\wedge^4V)}
\begin{document}

\title{ Addendum to "$K3$-surfaces of genus 8 and 
varieties of sums of powers of cubic fourfolds"}

\author{Atanas Iliev and Kristian Ranestad}

%
\date{ }
\maketitle

\vspace{1.0cm}

\begin{abstract}
In this note, which is an addendum to the e-print math.AG/9810121, we prove 
that the variety VSP(F,10) of presentations of a general cubic form F in 6 
variables as a sum of 10 cubes is a smooth symplectic 4-fold, which is 
deformation equivalent to the Hilbert square of a K3 surface of genus 8 but 
different from the family of lines on a cubic 4-fold. This provides a new 
geometric construction of a compact complex symplectic fourfold, different 
from a Hilbert square of a K3 surface, a generalized Kummer 4-fold, the 
variety of lines on a cubic 4-fold and the recent examples of O'Grady 
(see Duke Math. J. 134, no. 1 (2006), 99-137).  
\end{abstract}





Eyal Markman noticed that the proof of Theorem 3.17 in the paper \cite{IR} is 
incomplete. The purpose of this note is to correct both the 
statement and the proof of that theorem.  

In the paper we show that $Hilb_{2}S$ for a 
general $K3$ surface $S$ of genus $8$, coincides with the variety of 
sums of powers $VSP(F,10)$ of a cubic fourfold $F$.  On the other hand 
Beauville and Donagi showed that $Hilb_{2}S$ coincides with the 
variety of lines $Hilb_{l}(F')$ in a Pfaffian cubic fourfold $F'$.  The natural 
polarization on $VSP(F,10)$ and on $Hilb_{l}(F')$ provides 
two polarizations $L_{vsp}$ and $L_{lines}$ on $Hilb_{2}S$.  Now the 
pairs $(Hilb_{2}S,L_{vsp})$ and $(Hilb_{2}S,L_{line})$ both deform 
smoothly with the cubic $F$ (resp. $F'$), giving components in the 
moduli space of polarized hyperkahler manifolds.  The general element 
of each component has Picard group generated by the given polarization.
In the proof of Theorem 3.17  in \cite{IR} we overlooked the importance of the 
polarization when we concluded that the two components of the 
moduli space coincide.  In fact we prove 

\bigskip\noindent
{\bf Theorem.} The two polarizations $L_{vsp}$ and $L_{lines}$ are 
linearly independent on $Hilb_{2}S$.

\bigskip\noindent
By Huybrechts (\cite{H} 1.8)  and Lemma 3.18 of \cite{IR} we get the

\bigskip\noindent
{\bf Corollary.} The moduli spaces of $VSP(F,10)$ and of $Hilb_{l}(F)$ for general cubic 
fourfolds form two distinct components of the moduli space of polarized 
hyperkahler manifolds. They intersect each other transversally  along a divisor 
representing manifolds  $Hilb_{2}S$, for a $K3$-surface $S$ of 
genus $8$.

\bigskip\noindent
{\bf Proof}.  To show that the polarizations are independent we show that 
the degrees 
of a particular rational curve  with respect to these two polarizations are relatively 
prime.  The polarization $L_{lines}$ is primitive, so this 
suffices to conclude.

\medskip\noindent
The test curve $C_{p}$ is a ruling in the diagonal on $Hilb_{2}(S)$,
i.e. the pencil of length two subschemes supported at a point $p\in 
S$.  This pencil is of course defined by the pencil of lines through 
$p$ in the tangent plane to $S$ at $p$. 
We show that this curve has degree 5 with respect to $L_{lines}$ and 
degree 3 with respect to $L_{vsp}$.

\bigskip\noindent
  Recall the setting in \cite{IR}.  We fix a vector space  $V\cong {\bf 
  C}^6$, and consider the Grassmannians $\Gr 
2V\subset{\bf P}(\wedge^2V)$ and $\Gr 
4V$ in the natural dual space ${\bf P}(\wedge^4V)$. 
  The surface $S$ is a linear section $S=\Pn 8\cap \Gr 
2V$.  We denote the $\Pn 8$ that span $S$ by $\Pn 8(S)$.  The 
cubic fourfold $F'$ is the linear section defined by $\Pn 5(S)=\Pn 
8(S)^{\bot}$ of the Pfaffian cubic in ${\bf P}(\wedge^4V)$ (the 
secant variety of $\Gr 4V$).  Since $\Gr 4V\subset {\bf 
P}(\wedge^4V)$ is arithmetically Gorenstein and has codimension 6, the 
quotient of its homogeneous coordinate ring defined by  $\Pn 
5(S)$ is an Artinian Gorenstein ring.  The socle degree of this ring is 
$3$.  The dual 
socle generator is therefore a cubic form.  This form defines the socalled apolar cubic fourfold $F$ associated to $S$. 

There are three polarizations on $Hilb_{2}S$ that show up in this 
setting.
First, $Hilb_{2}S$ is parameterized by the secant (and tangent) lines to $S$ (we may 
assume that $S$ contains no lines), so there is a polarization $L_{sec}$ 
defined by the Pl\"ucker embedding of 
the variety of secant lines to $S$.  Secondly, there is a 
polarization $L_{lines}$ defined by the Pl\"ucker embedding of the 
variety of lines in $F'$. The third polarization $L_{vsp}$ defines 
 the vsp-map 
$$
f:Hilb_2 \ S \rightarrow VSP(F, 10) \subset {\bf G}(4,\wedge^2 \ V),
$$   
\noindent where $F$ is the apolar cubic 4-fold to $S$.  By \cite{IR} 
Lemmas 2.8 and 3.8, it is defined as follows: 
We first identify $l \in Hilb_2 \ S$ with a secant or tangent line to $S$. 
If $l$ is a secant line, then let $ \{ p,q \}=l \cap S $ and let 
$L_{p}$ and $L_{q}$ be the corresponding lines in ${\bf P}(V)$.
Since $S$ does not 
contain lines (in particular $S$ does not contain the line $<p,q>$ 
spanned by $p$ and $q$),  the span $<L_p \cup L_q> = {\bf 
P}(U_l)$ 
for some 4-space $U_l \subset V$. The set $Q_l$ of all lines
$L \subset {\bf P}(V)$ that intersect $L_p$ and $L_q$ 
is a smooth quadric surface $Q_l \subset {\bf G}(2,U_l)$, 
and the space ${\bf P}(l) = {\bf P}^3_l = Span \ Q_l \subset {\bf P}(\wedge^2 \ 
U_l)$ 
is just the image $f(l)$ of $l$. If $l$ is a tangent line, the same 
argument goes through for a quadric surface cone $Q_l$ with vertex 
at $p$ spanning a 
$3$-space ${\bf P}(l)$.
In both cases, note that ${\bf P}(l)$ is the polar in ${\bf P}(\wedge^2 \ 
U_l)$  of $l$ with respect 
to the quadric ${\bf G}(2,U_l)$. 

\bigskip\noindent
{\bf Lemma}.  $C_{p}\cdot L_{vsp}=3$.

\medskip\noindent
{\bf Proof}. In order to describe the restriction of the vsp-map $f$ 
to the curve 
$C_p\subset Hilb_2 \ S$, 
note first that the map 
$\beta : Hilb_2 \ S \rightarrow {\bf G}(4,V)$, 
$\beta : l \rightarrow U_l$ defined above, 
is the restriction on $Hilb_2 \ S \subset {\bf P}(\wedge^2 \ V)$
of the Cremona transformation 
$\beta: {\bf P}(\wedge^2 \ V) \rightarrow {\bf P}(\wedge^4 \ V)$ 
defined by the quadratic Pfaffians defining $\Gr 2V$, see \cite{IR}. 
The restriction $\beta_p$ of $\beta$ to the tangent plane 
${\bf P}^2_p$ to $S$ at $p$ is defined by the system of all conics
in ${\bf P}^2_p$ singular at $p$. Therefore $\beta_p$ blows up $p$ 
and maps the exceptional divisor $E_p$ over $p$ to a conic 
$q_p \subset {\bf G}(4,V)$. In the embedding 
$Hilb_2 \ S \subset {\bf P}(\wedge^2(\wedge^2 \ V))$ defined above,  
$C_p$ is identified with the pencil of lines through $p$ in the 
tangent plane ${\bf P}^2_p$, i.e. with the exceptional curve $E_{p}$.
Therefore the Cremona transformation 
$\beta$ sends the line $C_p \subset ({\bf P}^2_p)^*$ to the conic 
$q_p \subset {\bf G}(4,V)$ as above.

 The pencil 
of ${\bf P}(\wedge^2 \ U_l)$'s in ${\bf P}(\wedge^2 \ V)$ 
corresponding to $l\in C_{p}$
form a rational scroll $W$ of degree $6$ parameterized by the 
elements ${\beta}(l) \in q_p$ of the conic $q_p \subset {\bf G}(4,V)$ 
(the rational scroll of ${\bf P}(\wedge^2 \ U_l)$'s corresponding to a 
line in ${\bf G}(4,V)$ has degree $3$). 
The scroll $W$ has the point $p$ as its vertex, so $W$ contains 
the tangent plane ${\bf P}^2_p$ to $S$ at $p$, and for any $l$ 
the ruling ${\bf P}(\wedge^2 \ U_l) \subset W$ intersects ${\bf P}^2_p$
at the line $l$. 
Finally, for each $l\in C_{p}$ the space ${\bf P}(l) = {\bf P}^3_l = f(l)$ 
coincides with the polar $3$-space to the line 
$l \subset {\bf P}(\wedge^2 \ U_l)$ with respect to the quadric 
$\Gr 2{U_{l}}$, i.e. the isomorphic image 
$f(C_p) \subset VSP(F(S),10) \subset {\bf G}(4,\wedge^2 \ V)$
is the polar scroll to the plane ${\bf P}^2_p$ with respect to
the quadric bundle $W \cap {\bf G}(2,V)$. 
We use the following lemma to show that the degree of this polar 
scroll is $3$.

\bigskip\noindent
{\bf Lemma}  Let $X=\Pn n_{\Pn 1}$ be a $\Pn n$-bundle with $H$ the 
divisor associated to ${\cal O}_{X}(1)$ and $F$ the class of a fiber.  
Let $Q\subset X$ be a quadric bundle such that $Q\equiv 2H+aF$, and 
let $s:\Pn 1\to X$ be a section of the projective bundle such that 
$s^{*}({\cal O}_{X}(1))$ has degree $d$. Let $Y_{s}\subset X$ be
the variety in $X$ of polar $\Pn {n-1}$-spaces to $s(t)$ with respect 
to $Q_{t}$ where $t\in \Pn 1$.  Then $$Y_{s}\equiv H+(a+d)F.$$

\medskip\noindent
{\bf Proof}. This is straightforward computation.
Let $<x_{i}>$ be a basis for the sections of ${\cal O}_{X}(1)$, and let 
$u_{0},u_{1}$ generate $H^0(\Pn 1,{\cal O}_{\Pn 1}(1))$.  Then the equation 
for $Q$ is quadratic in the $x_{i}$ and of degree $a$ in the $u_{i}$.
The point $s(t)$ has coordinates of degree $d$ in the $u_{i}$.  By 
differentiation the polar $Y_{s}$ has a linear equation in the 
$x_{i}$ of degree $d+a$ in the $u_{i}$. \qed

\medskip\noindent
In our situation we first note that the quadric bundle $Q\subset W\subset G$ has 
degree $10$.  In fact the quadric bundle of quadrics $\Gr 2{U_{l}}$'s in the cubic scroll
of the ${\bf P}(\wedge^2 \ U_l)$'s corresponding to a 
line in ${\bf G}(4,V)$ is isomorphic to a tangent hyperplane section 
of a $\Gr 25$, and therefore has degree 5.  Since $Q$ is the 
corresponding quadric bundle in the sextic scroll $W$ corresponding to 
a conic in ${\bf G}(4,V)$, it has degree 10.
If $H$ is the 
pullback (to a desingularization of $W$) of the Pl\"ucker divisor 
on $\Gr 2V$, and if
$F$ is a member of the ruling, then 
$Q\equiv 2H-2F$ (on the desingularization of $W$).
Let $s_{0}:\Pn 1\to W$ be the constant map to the vertex $p$.  Then, by 
the lemma,
$$Y_{s_{0}}\equiv H-2F.$$
Let $Q_{0}=Q\cap Y_{s_{0}}\subset  Y_{s_{0}}=X_{0}$ and let $s_{1}:\Pn 
1\to X_{0}$ be the section corresponding to a line in the tangent plane 
that does not pass through the tangency point $p$.  Applying the 
lemma once more we get: 
$$Y_{s_{1}}\equiv H-F$$
on $X_{0}$.  Thus $Y_{s_{1}}$ has degree $$H^4\cdot (H-F)\cdot 
(H-2F)=3.$$
\qed

\bigskip\noindent
{\bf Lemma}.  $C_{p}\cdot L_{lines}=5$.

\medskip\noindent
{\bf Proof}. Let $V_{p}$ be the $2$-space corresponding to $p\in \Gr 2V$.  Each 
 tangent 
line $l$ 
to $\Gr 2V$ at $p$ is tangent to a unique quadric $\Gr 2{U_{l}}$ for a 
$4$-space $U_{l}$ that contains $V_{p}$.  The hyperplanes that contain 
$\Gr 2{U_{l}}$ form a $\Pn 8({U_{l}})$ in the dual space $\Pr$, and 
this $\Pn 8({U_{l}})$ is the tangent space 
to $\Gr 4V$ at the point corresponding to ${U_{l}}$.  The tangent space
$\Pn 8({U_{l}})$ is contained in the Pfaffian cubic of secant lines to 
$\Gr 4V$, and it intersects $\Pn 5(S)$ in a line $l'$, since $\Pn 5(S)^{\bot} 
=\Pn 8(S)$ intersects ${\bf P}(\wedge^2{U_{l}})=\Pn 8({U_{l}})^{\bot}$ in a line.
 We want to compute the degree of the curve of such lines $l'$ in $\Pr$ corresponding to the pencil of tangent lines 
to $S$ at $p$. 

Above we computed the degree of $W\subset\Pd$, the pencil of  ${\bf 
P}(\wedge^2U_{l})$'s that intersect $\Pn 8(S)$ in a tangent line to $S$.
The corresponding pencil of $\Pn 8(U_{l})$'s orthogonal to the ${\bf 
P}(\wedge^2U_{l})$'s  is a  rational  $9$-fold scroll $W'$ of degree $6$ in 
$\Pr$. Since $W$ is a cone, the scroll $W'$ is contained in 
a hyperplane.  Furthermore $\Pn 5(S)$ is contained in this 
hyperplane.  In fact, since $W\cap \Pn 8(S)$ is precisely a plane (the 
tangent plane to $S$)
through the vertex $p$ that intersect each ruling of $W$ in a line, 
the intersection of each $\Pn 8(U_{l})$ with $\Pn 5(S)$ is a line.
The degree of $W'\cap \Pn 5(S)$ coincides with the degree of the 
projection of $W$ from the tangent plane.
But this degree is $5$, so the lemma follows.
\qed

\noindent {\bf Acknowledgement.} We thank Brendan Hassett for helpful discussions on the Fano variety 
of lines in a cubic fourfold.

\end{document}